\documentclass[11pt]{amsart}

\usepackage{amsmath,amssymb,amsthm,amsrefs,a4wide}
\usepackage{latexsym}
\usepackage{indentfirst}
\usepackage{graphicx}
\usepackage{placeins}
\usepackage{booktabs}
\usepackage{algorithm}
\usepackage{algorithmic}
\usepackage{multirow}
\usepackage{color}

\theoremstyle{plain}

\theoremstyle{definition}

\newtheorem{example}{Example}

\theoremstyle{remark}

\newcommand{\diag}{\mathrm{diag}}

\newcommand{\bbm}{\begin{bmatrix}}
\newcommand{\ebm}{\end{bmatrix}}

\newcommand{\C}{\mathbb{C}}

\newcommand{\x}{\tilde{x}}
\newcommand{\w}{\tilde{w}}

\begin{document}

\title[Perturbative analysis for noisy spectral estimation]{A perturbative analysis for noisy spectral estimation}

\author[]{Lexing Ying} \address[Lexing Ying]{Department of Mathematics, Stanford University,
  Stanford, CA 94305} \email{lexing@stanford.edu}

\thanks{This work is partially supported by NSF grant DMS-2011699 and DMS-2208163.}

\keywords{Perturbative analysis, spectral estimation.}

\subjclass[2010]{65R32,65T99}

\begin{abstract}
Spectral estimation is a fundamental task in signal processing. Recent algorithms in quantum phase estimation are concerned with the large noise, large frequency regime of the spectral estimation problem. The recent work in Ding-Epperly-Lin-Zhang shows that the ESPRIT algorithm exhibits superconvergence behavior for the spike locations in terms of the maximum frequency. This note provides a perturbative analysis to explain this behavior. It also extends the discussion to the case where the noise grows with the sampling frequency.
\end{abstract}

\maketitle

%----------------------------------------------------------
\section{Introduction}\label{sec:intro}

%---------
This note is concerned with the spectral estimation problem, a fundamental task in signal processing.  Let us consider an unknown spectrum measure of the form
\[
\mu(x) \equiv \sum_{k=1}^r \delta_{x_k}(x) w_k,
\]
where $\{x_k\}\subset [0,2\pi)$ are the spike locations and $\{w_k\}$ are the weights. Let $f_j$ be the Fourier transform of the measure $\mu$ at integer $j$, i.e.,
\[
f_j = \sum_{k=1}^r e^{i j x_k} w_k.
\]
Suppose that we have access to the noisy measurement $g_j \equiv f_j + z_j$ for $-n \le j \le n$, where $z_j$ is the measurement noise. In what follows, we assume that each $z_j$ is an independent complex Gaussian random variable $\sigma N_\C(0,1)$ with standard deviation $\sigma$, denoted as $z_j \sim \sigma N_\C(0,1)$. The spectrum estimation problem is to estimate the locations $\{x_k\}_{1\le k \le r}$ and the weights $\{w_k\}_{1\le k \le r}$ based on the noisy $\{g_j\}_{-n \le j \le n}$.

Spectral estimation is a fundamental problem in signal processing. Over the decades, many theoretical studies have been devoted to this problem. When the noise level is relatively small, this is often referred to as the superresolution setting \cite{candes2014towards,donoho1992superresolution,demanet2015recoverability,moitra2015super,li2020super,liao2016music}. The algorithms proposed for this include optimization-based approaches and algebraic methods such as \cite{prony1795essai,roy1989esprit,hua1990matrix}

In recent years, spectral estimation has become an important tool in quantum computing, especially in quantum phase estimation \cite{ding2023even,ding2023simultaneous,li2023adaptive,li2023note,ni2023low,ding2024quantum}. However, the setup is somewhat different: the noise level is relatively large, but it can be compensated by collecting Fourier measurements for large values of $j$. Therefore, this is the large noise, large $n$ regime.

Recently in \cite{ding2024esprit}, Ding et al showed that, for mean-zero $O(1)$ noise, the ESPRIT algorithm \cite{roy1989esprit} results in $O(n^{-3/2})$ error in $\{x_k\}$ and $O(n^{-1/2})$ error in $\{w_k\}$. The superconvergence result $O(n^{-3/2})$ for $\{x_k\}$ is quite surprising and the proof in \cite{ding2024esprit} is a tour de force.

%{\bf Contributions.}

The goal of this short note is to provide an intuitive understanding of the $O(n^{-3/2})$ error scaling in $\{x_k\}$ and $O(n^{-1/2})$ error scaling in $\{w_k\}$. We conduct a perturbative analysis to explain this superconvergence behavior of $\{x_k\}$. This argument is also extended to the more general case where the standard deviation of the noise grows with frequency, which is highly relevant for quantum phase estimation.

The rest of the paper is organized as follows. Section \ref{sec:pa} presents the perturbative error analysis for the large noise large $n$ setting. Section \ref{sec:ge} considers the more general noise scaling model. Numerical results are also provided for different noise scalings to compare with the perturbative analysis.

%----------------------------------------------------------
\section{Perturbative analysis} \label{sec:pa}

The analysis in this section makes two assumptions.
\begin{itemize}
\item The minimum gap between $\{x_k\}$ is bounded from below by a constant $\Delta>0$.
\item Each $w_k$ is $\Theta(1)$.
\end{itemize}
From the perspective of maximum likelihood estimation (MLE), one approach for spectral estimation is to solve for the approximations $\{\x_k\}$ and $\{\w_k\}$ from the following optimization problem
\[
\min_{\x_k,\w_k} \sum_{j=-n}^n \left|\sum_{k=1}^r e^{ij \x_k} \w_k - g_j \right|^2. 
\]
Let us introduce
\[
\x_k = x_k + a_k, \quad \w_k = w_k (1+b_k)
\]
where $a_k$ is the absolute error of $x_k$ and $b_k$ is the relative error of $w_k$. In terms of $a_k$ and $b_k$, the residue for each $j$ is
\[
\sum_{k=1}^r e^{ij \x_k} \w_k - g_j = \sum_{k=1}^r e^{ij x_k} w_k e^{ija_k} (1+b_k) - (f_j+z_j).
\]
Expanding in terms of $\{a_k\}$ and $\{b_k\}$ and ignoring the higher order terms gives
\[
\sum_{k=1}^r e^{ij x_k} w_k (1+ ija_k + b_k)  - (f_j+z_j).
\]
Using the idenitity $\sum_{k=1}^r e^{ij x_k} w_k = f_j$ for exact $\{x_k\}$ and $\{w_k\}$ leads to the leading order difference
\[
\sum_{k=1}^r e^{ij x_k} (ij a_k + b_k) - z_j.
\]
Therefore, modulus the higher order terms, the optimization problem becomes a least square problem in terms of $\{a_k\}$ and $\{b_k\}$
\[
\min_{a_k,b_k} \sum_{j=-n}^n \left|\sum_{k=1}^r e^{ij x_k} w_k (ij a_k + b_k) - z_j \right|^2. 
\]
In order to put this into a matrix form, let us define the matrices
\[
A=\bbm ij e^{ij x_k} \ebm_{-n \le j\le n,1\le k\le r}, \quad B = [e^{ij x_k}]_{-n \le j\le n,1\le k\le r}, \quad W = \diag([w_k]_{1\le k\le r}),
\]
and 
\[
M =\bbm A & B \ebm \bbm W &  \\ & W\ebm.
\]
Define the vectors $z=\bbm z_j\ebm_{-n \le j \le n}$, $a=\bbm a_k \ebm_{1\le k\le r}$ and $b=\bbm b_k \ebm_{1\le k\le r}$. Then, the linear system from the least square problem is 
\[
M^* M \bbm a \\b\ebm = M^* z.
\]
and, from the factorization of $M$, 
\begin{equation}
  \bbm A^*A & A^* B \\ B^* A & B^* B\ebm  \bbm W &  \\ & W\ebm  \bbm a \\b\ebm  = \bbm A^* \\ B^*\ebm z
  \label{eq:lseq}
\end{equation}

Let us consider the matrix
\[
\bbm A^*A & A^* B \\ B^* A & B^* B\ebm.
\]
\begin{itemize}
\item For the top-left block $A^* A$, the $(k',k)$-th entry is $\sum_{j=-n}^n j^2 e^{ij (x_k-x_k')}$. Each diagonal entry is $\frac{2}{3} n^3 (1+O(\frac{1}{n}))$. To estimate the off-diagonal ones, we denote $t=x_k-x_k'$, consider
\begin{equation}
  \sum_{j=-n}^n e^{ij t} = \frac{\sin(n+\frac{1}{2})t}{\sin(\frac{1}{2})t},
  \label{eq:sum}
\end{equation}
and take the second-order derivative in $t$. Bounding the derivative shows that each off-diagonal entry is $O(\frac{n^2}{\Delta} + \frac{1}{\Delta^3})$. Since we care about the large $n$ behavior, assuming $n\gtrsim \frac{1}{\Delta}$ ensures that every off-diagonal entry is bounded by $O(\frac{n^2}{\Delta})$.

\item For the top-right block $A^* B$, the $(k',k)$-th entry is $\sum_{j=-n}^n (-ij) e^{ij (x_k-x_k')}$. Each diagonal entry is zero due to cancellation, while each off-diagonal one is $O(\frac{n}{\Delta})$ by considering the first order derivative of \eqref{eq:sum}.

\item The bottom-left block $B^* A$ is simillar to $A^* B$.

\item For the bottom-right block $B^* B$, the $(k',k)$-th entry is $\sum_j e^{ij (x_k-x_k')}$. Each diagonal entry is $2n (1+O(\frac{1}{n}))$, while each off-diagonal entry is $O(\frac{1}{\Delta})$.
\end{itemize}
Therefore, the whole matrix takes the form
\[
\bbm A^*A & A^* B \\ B^* A & B^* B\ebm =
\bbm \frac{2}{3} n^3 (1+O(\frac{1}{n})) I + O(n^2/\Delta) & O(n/\Delta) \\ O(n/\Delta) & 2n (1+O(\frac{1}{n})) I + O(1/\Delta) \ebm.
\]
When $n \gtrsim \frac{r}{\Delta}$, the diagonal entries dominate the non-diagonal ones, leading to the factorization
\[
\bbm A^*A & A^* B \\ B^* A & B^* B\ebm =
\bbm n^{3/2}I & \\ & n^{1/2} I \ebm
\left( \bbm \frac{2}{3} I & \\ & 2I \ebm + o(1) \right)
\bbm n^{3/2}I & \\ & n^{1/2} I \ebm
\]
and its inverse
\begin{equation}
  \bbm A^*A & A^* B \\ B^* A & B^* B\ebm^{-1} =
  \bbm n^{-3/2}I & \\ & n^{-1/2} I \ebm
  \left( \bbm \frac{3}{2} I & \\ & \frac{1}{2}I \ebm + o(1) \right)
  \bbm n^{-3/2}I & \\ & n^{-1/2} I \ebm.
  \label{eq:invfac}
\end{equation}

Let us next consider $\bbm A^* \\ B^*\ebm z = \bbm A^* z\\ B^* z\ebm$.
\begin{itemize}
\item The $k$-th entry of $A^*z$ is $\sum_j (-ij) e^{-ij x_k} z_j$. By the independence of $z_j \sim \sigma N_\C(0,1)$, each entry of the vector $A^*z$ is $\sim \sqrt{\frac{2}{3}} n^{3/2} \sigma N_\C(0,1)$.
\item The $k$-th entry of $B^*z$ is $\sum_j e^{-ij x_k} z_j$. The same argument shows that each entry of the vector $B^*z$ is $\sim \sqrt{2} n^{1/2} \sigma N_\C(0,1)$.
\end{itemize}

To solve \eqref{eq:lseq}, we introduce
\[
\bbm u \\ v \ebm
\equiv
\bbm A^*A & A^* B \\ B^* A & B^* B\ebm^{-1} \bbm A^* z \\ B^* z \ebm.
\]
Using the factorization \eqref{eq:invfac}, each entry of $u$ is $\sim \sqrt{\frac{3}{2}} n^{-3/2} \sigma N_\C(0,1)$, while each entry of $v$ is $\sim \sqrt{\frac{1}{2}} n^{-1/2} \sigma N_\C(0,1)$. Finally, because
\[
\bbm a \\b\ebm = 
\bbm W^{-1} &  \\ & W^{-1} \ebm  \bbm u \\ v \ebm
\]
and the weights $w_k$ are all $\Theta(1)$, we conclude that each entry of $a$ is $\sim O(n^{-3/2}) \sigma N_\C(0,1)$, while each entry of $b$ is $\sim O(n^{-1/2}) \sigma N_\C(0,1)$.

To summarize, when focusing the $n$ dependence and considering the regime $n \gtrsim \frac{r}{\Delta}$, the absolute errors of the spike locations $\{x_k\}$ scale like $O(n^{-3/2})$, while the relative errors of the weights $\{w_k\}$ scale like $O(n^{-1/2})$.

%----------------------------------------------------------
\section{General noise} \label{sec:ge}

We have so far assumed that the noise $z_j \sim \sigma N_\C(0,1)$ has a variance independent of $j$. In the quantum phase estimation problems, a noise magnitude can grow with $j$ since it takes a longer quantum circuit to prepare. Therefore, a more realistic noise model is
\begin{equation}
  z_j \sim j^{p} \cdot \sigma N_\C(0,1)
  \label{eq:nnm}
\end{equation}
for some $p>0$. The above analysis can be easily adapted to this more general noise model. The only difference concerns the right hand side $ \bbm A^* \\ B^*\ebm z = \bbm A^* z \\ B^* z\ebm$. With the model \eqref{eq:nnm}, each entry of $A^* z$ has a standard deviation $O(n^{3/2+p}\sigma)$, while the one of $B^* z$ has a standard deviation $O(n^{1/2+p}\sigma)$. Therefore, for
\[
\bbm a \\b\ebm = 
\bbm W^{-1} &  \\ & W^{-1} \ebm
\bbm A^*A & A^* B \\ B^* A & B^* B\ebm^{-1} \bbm A^* z \\ B^* z \ebm,
\]
each entry of $a$ is $\sim O( n^{-3/2+p} \sigma) \cdot N_\C(0,1)$, while each one of $b$ is $\sim O(n^{-1/2+p} \sigma) \cdot N_\C(0,1)$.  This states that, as long as the $p<1/2$, the errors of $\{w_k\}$ decay as $n$ grows.  As long as $p<3/2$, the errors of $\{x_k\}$ goes to zero with $n$.

Below, we provide several numerical examples with different noise scalings. In each case, $r=4$, $\Delta=0.1$, $\sigma=0.1$, and multiple $n$ values are tested, ranging from $16$ to $4096$.

\begin{example}
  $p=0$, and this is the standard setting of Section \ref{sec:pa}. Figure \ref{fig:EX1} gives the error asymptotics for $\{x_k\}$ and $\{w_k\}$ as $n$ varies. The slopes of the log-log plots match well with the predictions $-1.5$ and $-0.5$ of the perturbative analysis, respectively.
  \begin{figure}[h!]
    \centering
    \includegraphics[scale=0.35]{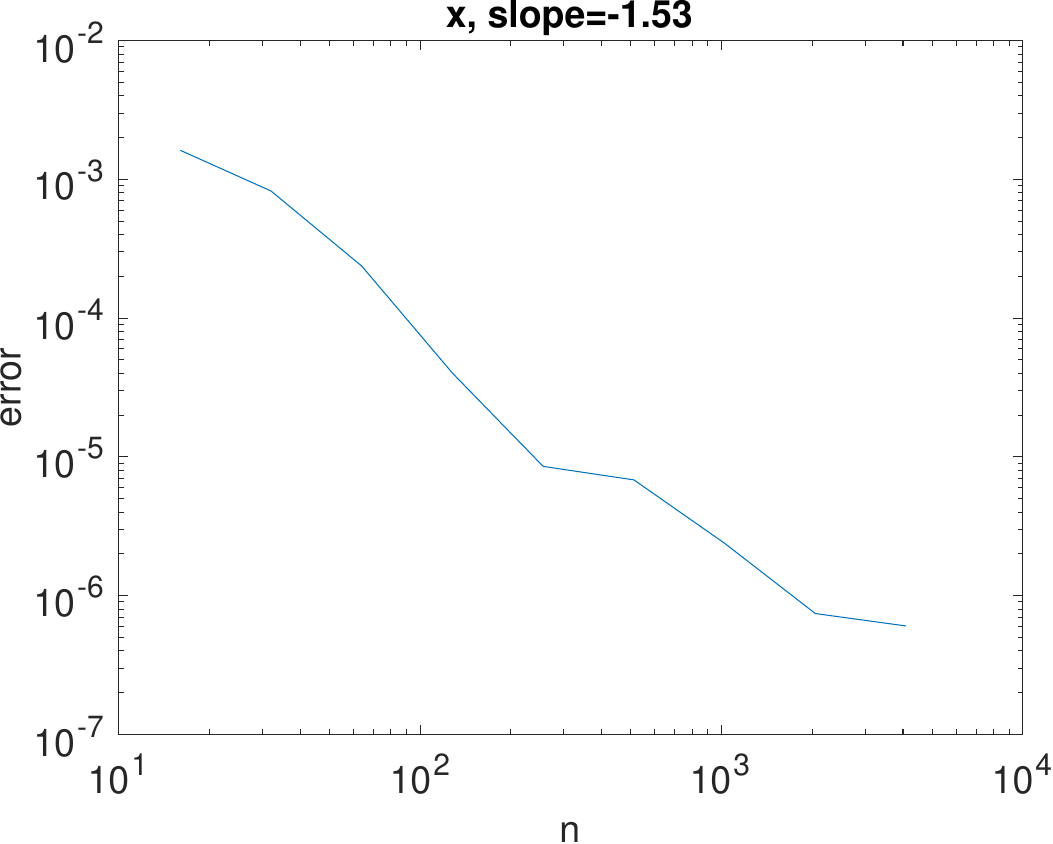}
    \includegraphics[scale=0.35]{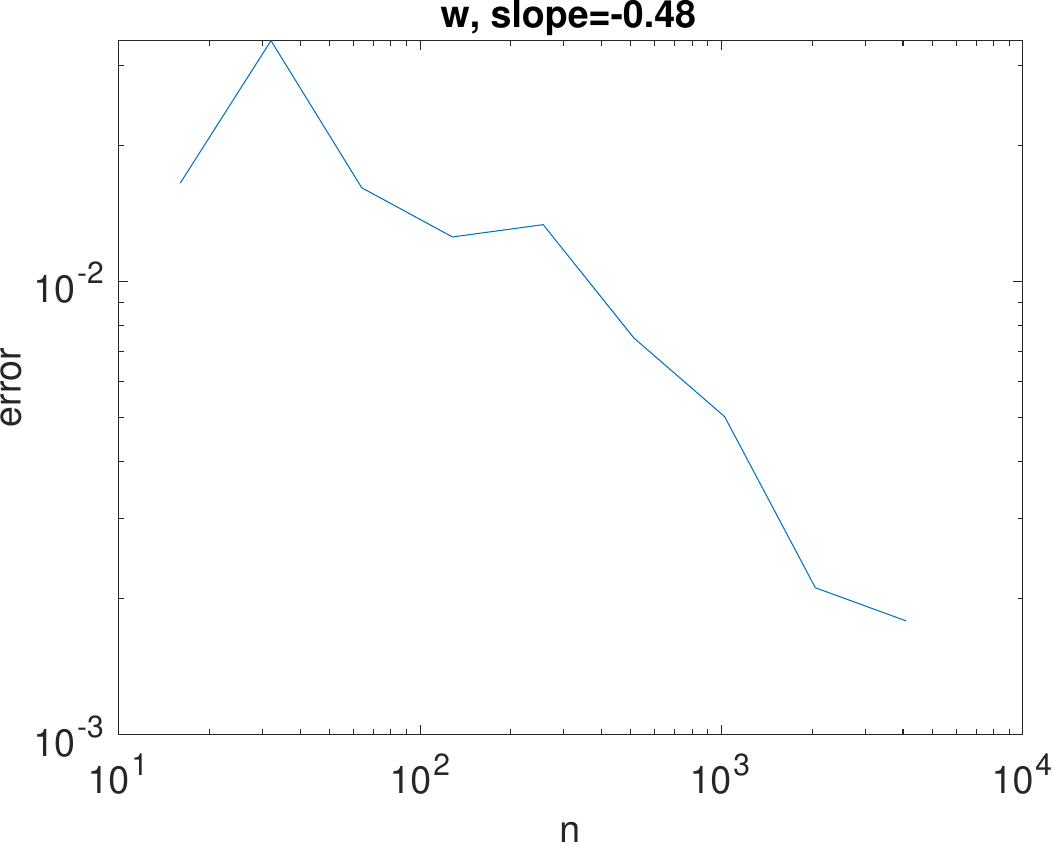}
    \caption{$p=0$. The slopes of the log-log plots for $\{x_k\}$ and $\{w_k\}$ match the perturbative analysis.}
    \label{fig:EX1}
  \end{figure}
\end{example}

\begin{example}
  $p=0.25$. Figure \ref{fig:EX2} gives the error asymptotics for $\{x_k\}$ and $\{w_k\}$ as $n$ varies. The slopes of the log-log plots also match well with the analytical predictions $-1.25$ and $-0.25$, respectively. Though the noise grows with $j$, the errors in $\{x_k\}$ and $\{w_k\}$ still exhibit the predicted decay rates.
  \begin{figure}[h!]
    \centering
    \includegraphics[scale=0.35]{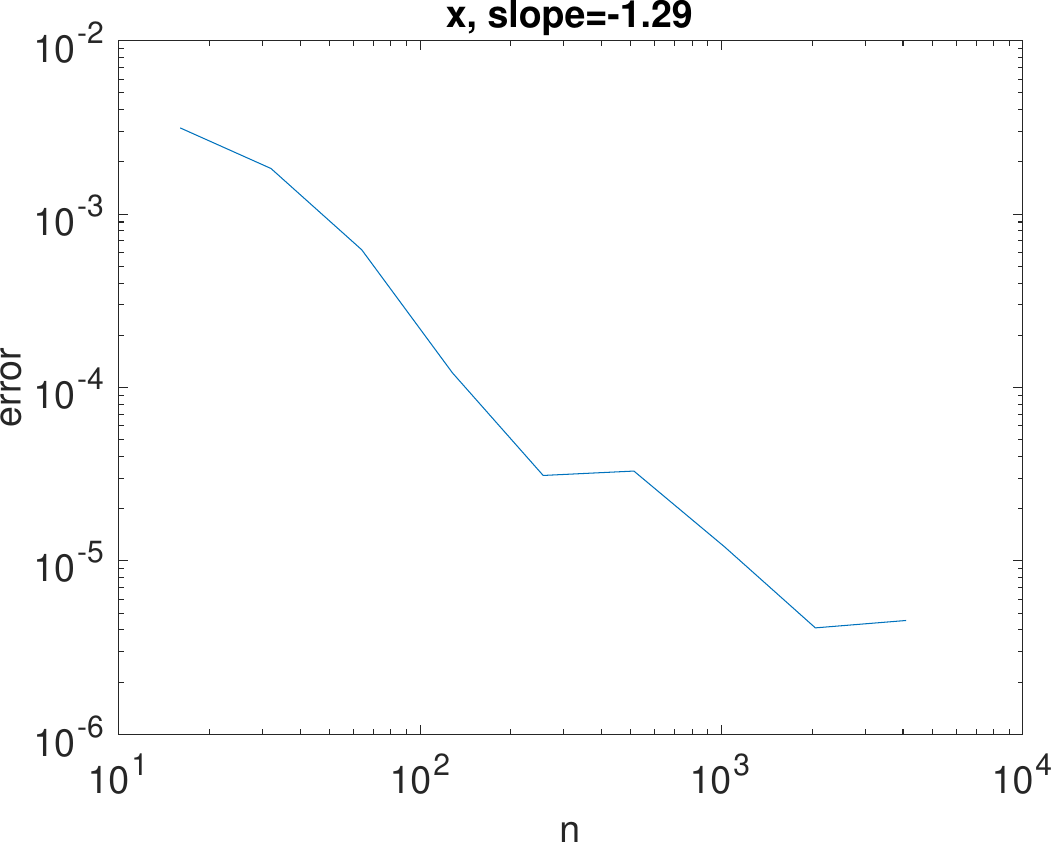}
    \includegraphics[scale=0.35]{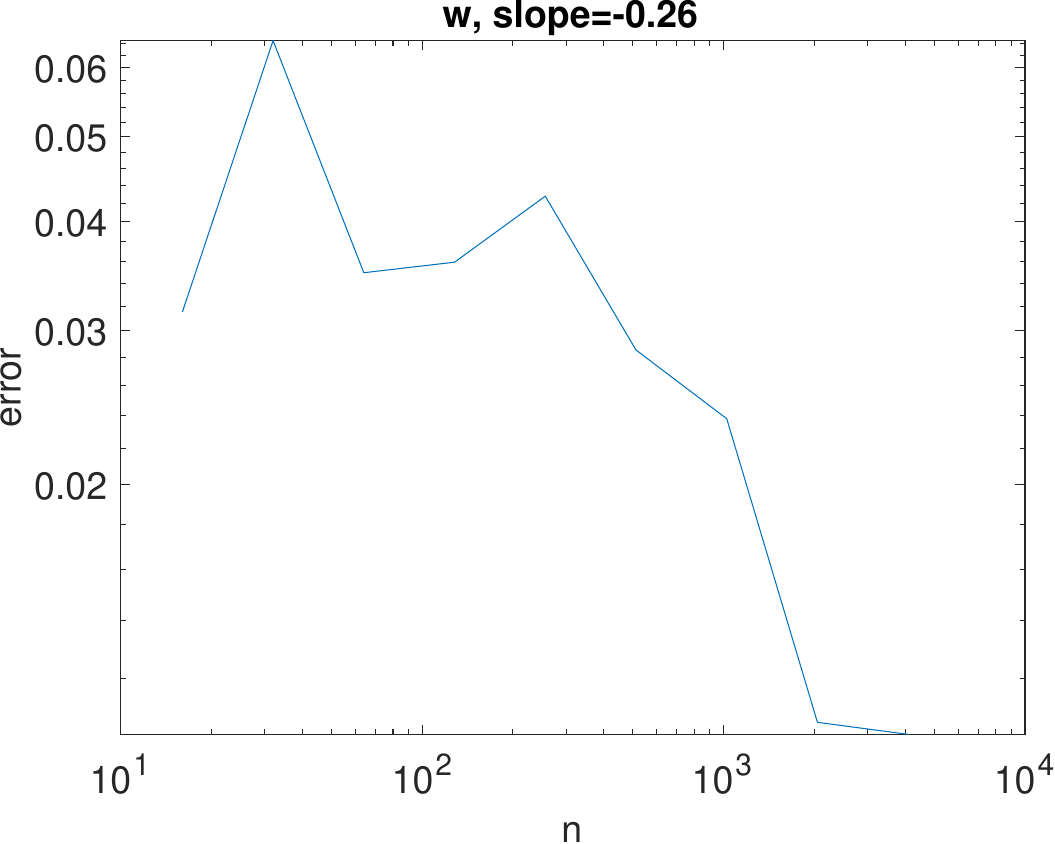}
    \caption{$p=0.25$. The slopes of the log-log plots for $\{x_k\}$ and $\{w_k\}$ match the perturbative analysis.}
    \label{fig:EX2}
  \end{figure}
\end{example}

\begin{example}
  $p=0.75$. Since $p>0.5$, we expect the error of $\{x_k\}$ to decay while the one of $\{w_k\}$ to grow. Figure \ref{fig:EX3} summarizes the error asymptotics for $\{x_k\}$ and $\{w_k\}$ as $n$ varies. The slopes of the log-log plots again match well with the prediction
  $-0.75$ and $0.25$, respectively.
  \begin{figure}[h!]
    \centering
    \includegraphics[scale=0.35]{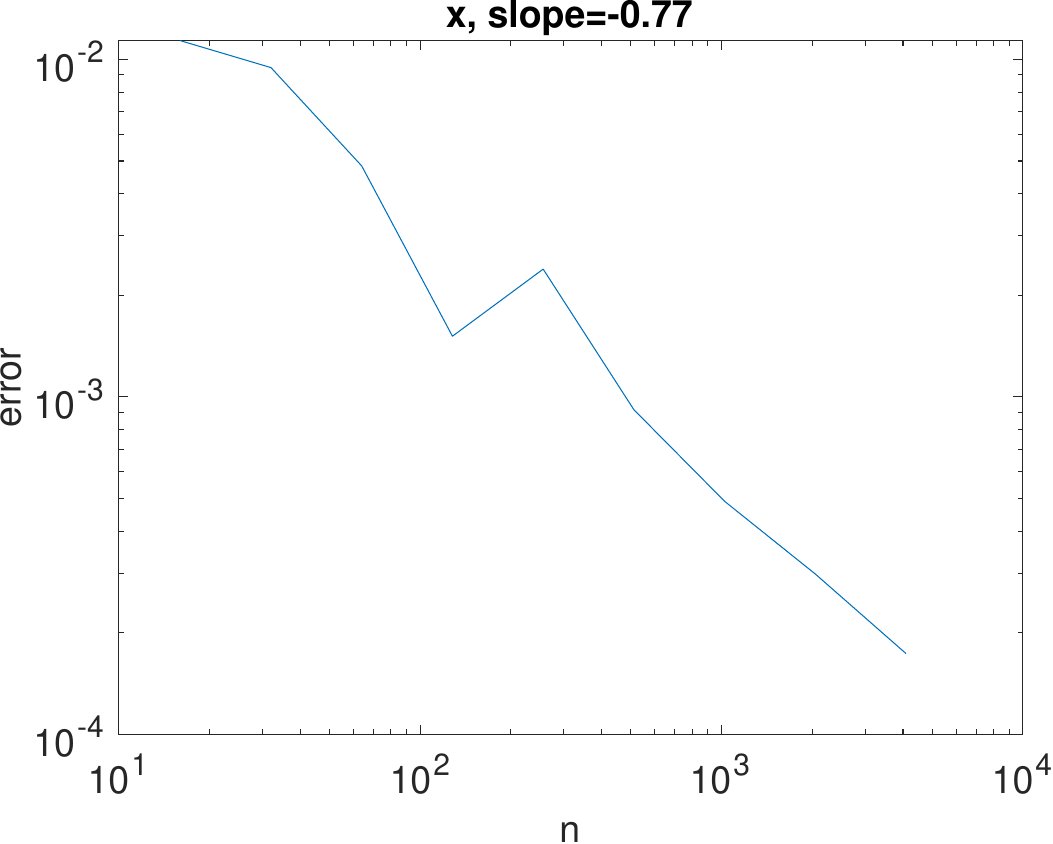}
    \includegraphics[scale=0.35]{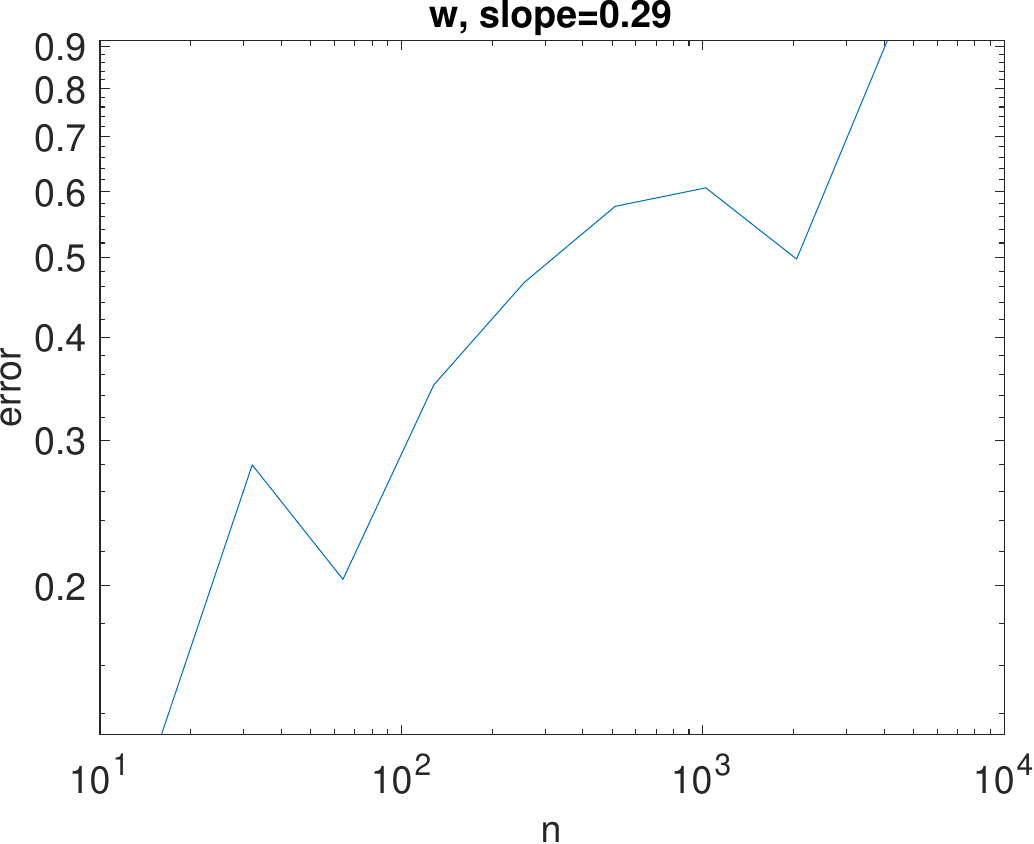}
    \caption{$p=0.75$. The slopes of the log-log plots for $\{x_k\}$ and $\{w_k\}$ match the perturbative analysis.}
    \label{fig:EX3}
  \end{figure}
\end{example}

\bibliographystyle{abbrv}

\bibliography{ref}

\end{document}